\documentclass[12pt]{amsart}
\usepackage{amssymb}

  \def\k{\Bbbk} \def\B{\mathcal B}
   
\def\tr{\triangle}

\def\End{\operatorname{Lin}} 
 
 \newcommand{\id}{\operatorname{id}} 
  
 \def\lin{\operatorname{lin}}
 \def\Mkg{\mathbf{bialg}}
\def\Mkg{{\mathcal M}^{kG}}
\def\B{\mathcal B}
\def\bbox{\qquad\square} \def\prof{\noindent{\em Proof. }}

\begin{document}

\numberwithin{equation}{section} \newtheorem{defn}{Definition}[section]
\newtheorem{lem}[defn]{Lemma} \newtheorem{prop}[defn]{Proposition}
\newtheorem{thm}[defn]{Theorem}

\title[Bicharacters, braids and Jacobi identity] {Bicharacters, braids and
Jacobi identity} \author{Jerzy R\'o\.za\'nski} \thanks {The work was
supported in part by KBN grant \# 2-P302-023-07} 
\address{Department of Mathematical Methods in Physics, University of Warsaw,
Ho\.za 74, PL-00-682 Warszawa, Poland.} \email{rozanski@fuw.edu.pl}

\begin{abstract} For an abelian group $G$ we consider braiding in a category
of $G$-graded modules $\Mkg$ given by a bicharacter $\chi$ on $G$. For
$(G,\chi)$-bialgebra $A$ in $\Mkg$ an analog of Lie bracket is defined. This
bracket is determined by a linear map $E\in\End(A)$ and n-ary operations
$\Omega^{n}_{E}$ on $A$.  Our result states that if $E(1)=0,\;E^{2}=0$ and
$\Omega^{3}_{E}=0$ then a braided Jacobi identity holds and the linear map
$E$ is a braided derivation of a braided Lie algebra.  \end{abstract}
\maketitle

\section{Introduction}

For an abelian group $G$ a monoidal category of $G$-graded modules $\Mkg$ has
braiding $\B$ given by a bicharacter $\chi$ on the group $G$.  In \cite{BP}
Pareigis considered a theory of n-ary Lie algebras in a braided monoidal
category for a model of braiding which is generated by a bicharacter. For
some bicharacters Pareigis defined a n-ary bracket with two types of Jacobi
identities.  For a binary operations a braided Lie algebra was considered as
a generalization of formal Lie theory by Gurevich for symmetrical braiding
\cite{gurevich} and for pure braiding a braided Lie bialgebra was introduced
by Majid \cite{majid}.

The aim of this short letter is to study a relationship between braiding and
Jacobi identity.  Our viewpoint consists a characterization of this problem
by some differential condition due to Koszul \cite{Koszul}.

Consider $(G,\chi)$-bialgebras \cite{BP} in the category $\Mkg$ which are
$\B$-commutative.  We start from the definition of n-ary operation
$\Omega^{n}_{E}$ on a bialgebra for an arbitrary linear map $E$ but restrict
ourselves to some brackets defined by binary operations $\Omega^{2}_{E}$
only.  If $E(1)=0,\;E^{2}=0$ and $\Omega^{3}_{E}=0$ then the Jacobi identity
for the binary operation holds. In this case the linear map $E$ is a braided
derivation of a Lie bracket.  Moreover we construct a braided derivation of
an underlying algebra.

\section{Notations, Definitions and remarks}

Let $\k$ be a field and let $\k^{\ast}$ be the multiplicative group of $\k$.
Denote by $G$ an abelian group written additively, $g_{i}\in G,\;
i=1,\;2,\;3.$ \begin{defn}[\cite{rybki}]\label{bich} A function $\chi:G\times
G\rightarrow\k^{\ast},\; \chi_{g_{1},\;g_{2}}\in\k^{\ast} $ is a bicharacter
of $G$ if \begin{gather*}
\chi_{g_{1}+g_{2},\;g_{3}}=\chi_{g_{1},\;g_{3}}\chi_{g_{2},\;g_{3}},\qquad
\chi_{g_{1},\;g_{2}+g_{3}}=\chi_{g_{1},\;g_{2}}\chi_{g_{1},\;g_{3}}
\end{gather*} \end{defn} A bicharacter is symmetric when
$\chi_{g_{1},\;g_{2}}=\chi_{g_{2},\;g_{1}}^{-1}.$ In this paper we will use
nonsymmetrical bicharacters.

For objects $L,\;M\in\Mkg$ we have the tensor product $$(L\otimes
M)_{g}={\oplus}_{h\in G} L_{h}\otimes M_{g-h}$$ and we have a braiding $\B$
for homogeneous elements $l\in L,\;m\in M$ with the corresponding degrees
$|l|,\;|m|$, \begin{gather}\label{braid}\B:L\otimes M\ni l\otimes
m\rightarrow \chi_{|l|,\;|m|}\;  m\otimes l\in M\otimes L.\end{gather} Note
that $\Mkg$ with $\B$ is a strict braided monoidal category \cite{rybki}.

Consider $(G,\;\chi)$-algebras in the category $\Mkg$. An algebra is a pair
$A=\{L,\;m\}$ where $L\in\Mkg$ and $m:L^{\otimes 2}\rightarrow L$, the
multiplication in the algebra $A$ preserves inclusions $L_{g}L_{h}\subset
L_{g+h}.$ For algebras $A,\;B$ the object $A\otimes B$ is again an algebra
with a multiplication as composition $$A\otimes B\otimes A\otimes
B\stackrel{\B}{\rightarrow} A\otimes A\otimes B\otimes
B\stackrel{m_{A}\otimes m_{B}}{\rightarrow} A\otimes B.$$ Above definition
works similarly with a (co)multiplication for categories of
$(G,\chi)$-coalgebras, bialgebras and Hopf algebras.

In \cite{OPR} we consider a braided derivations of an algebra.  In this paper
we have an interesting example of such derivations where a nonsymmetrical
braiding is constructed by a bicharacter on a group.
\begin{defn}[\cite{OPR}]\label{ader} A braided derivation of an associative
algebra $A=\{L,\;m\}$ is $\k$-linear map $d\in\lin_{\k}(L,\;L)$ that
satisfies the Leibniz rule $$d\circ m:=m\circ (d\otimes\id_{L})+m\circ
B^{-1}\circ (d\otimes\id_{L})\circ B$$ \end{defn} For consistent conditions
of derivations with the associativity of an algebra $A$ see in \cite{OPR}.

\section{Main Propositions}

Let $L\in\Mkg$ and $(G,\chi)$-bialgebra $A=\{L,\;m,\;\tr\}$ be
$\B$-commutative.  Let us fix the comultiplication $\tr:A\rightarrow
A^{\otimes 2}$ such that $\tr(a)=a\otimes 1-1\otimes a$.  The extension of
the comultiplication $\tr$ is a map $\tr^{n}:\otimes^{n} A\rightarrow
A\otimes A$ \begin{gather*} \tr^{n}(a_{1}\otimes...\otimes
a_{n})=\tr(a_{1})\cdot...\cdot\tr(a_{n}).  \end{gather*} For example
$\tr(a)\tr(b)=ab\otimes\id-a\otimes b -\chi_{a,\;b}b\otimes a+\id\otimes ab.$

For a linear map $E\in\End(A)$ one can define \cite{Koszul} the n-form
$\Omega^{n}_{E}:  A^{\otimes n}\rightarrow A$,
\begin{gather}\Omega_{E}^{n}(a_{1}\otimes...\otimes
a_{n})=m\circ(E\otimes\id)\tr^{n}(a_{1}\otimes... a_{n})..  \end{gather} For
example \begin{gather} \Omega_{E}^{2}(a\otimes
b)=E(ab)-E(a)b-\chi_{a,\;b}E(b)a+E(1)ab,\label{Om2}\\ \Omega_{E}^{3}(a\otimes
b\otimes c)=E(abc)-\chi_{b,\;c}E(ac)b-\chi_{a,\;b+c}E(bc)a\label{Om3}\\
+\chi_{a+b,\;c}E(c)ab-E(ab)c+E(a)bc+\chi_{a,\;b}E(b)ac-E(1)abc.\nonumber
\end{gather}

Let us assume that a bicharacter satisfies \begin{equation}\label{OmegaE2}
\chi_{a,\;b}\chi_{E(a),\;E(b)}\chi^{-1}_{a,\;E(b)}\chi^{-1}_{E(a),\;b}=-1.
\end{equation} \begin{lem}\label{Der1} For bicharacters (\ref{OmegaE2}) and
$E(1)=0$ we have $$\Omega^{2}_{E^{2}}(a\otimes b)=E\Omega^{2}_{E}(a\otimes b)
+\Omega^{2}_{E}(E(a)\otimes b)+
\chi_{a,\;b}{\chi}^{-1}_{a,\;E(b)}\Omega^{2}_{E}(a\otimes E(b)).  $$
\end{lem} \prof For the right side (\ref{Om2}) with (\ref{OmegaE2}) is
used.$\bbox$\\

The three-ary operations $\Omega^{3}_{E}$ can be expressed by the binary
$\Omega^{2}_{E}$.  \begin{lem}\label{der2}
$$\Omega^{3}_{E}(a\otimes b\otimes c)=\Omega^{2}_{E}(a\otimes
bc)-\Omega^{2}_{E}(a\otimes b)c-\chi_{b,\;c}\Omega^{2}_{E}(a\otimes c)b.$$
\end{lem} \prof From (\ref{Om3}) the right side has three terms
\begin{multline*}
\Omega^{2}_{E}(a\otimes bc)=E(abc)-E(a)bc-\chi_{a,\;bc}E(bc)a+E(1)abc,\\
-\Omega^{2}_{E}(a\otimes b)c=-E(ab)c+E(a)bc+\chi_{a,\;b}E(b)ac-E(1)abc,\\
-\chi_{b,\;c}\Omega^{2}_{E}(a\otimes c)b
=-\chi_{b,\;c}E(ac)b+\chi_{b,\;c}E(a)cb\\+\chi_{b,\;c}\chi_{a,\;c}E(c)ab
-\chi_{b,\;c}E(1)acb.
\end{multline*}
Due to the $\B$-commutativity the sum is $\Omega^{3}_{E}(a\otimes b\otimes
c)$.$\bbox$\\

Let us assume that a bicharacter satisfies \begin{gather}
\chi_{a,\;b+E(bc)}\chi^{-1}_{a,\;E(b)+b+c}=1,\quad
\chi_{b,\;E(c)}\chi_{a,\;E(bc)}\chi^{-1}_{a+b,\;E(c)}=1,\label{Jacobi}\\
\chi_{a+b,\;c}\chi_{\Omega(a,b),\;E(c)}\chi^{-1}_{a+b,\;E(c)}=-1,\quad
\chi_{a,\;b}\chi_{E(b),\;c}\chi_{\Omega^{2}_{E}(a,c),\;
E(b)}\chi^{-1}_{a,\;E(b)}\chi^{-1}_{b,\;c}=-1\nonumber.
\end{gather}

\begin{lem}\label{Jacobi1} For bicharacters (\ref{OmegaE2}), (\ref{Jacobi})
and $E(1)=0$ we have
\begin{multline*} \Omega^{2}_{E}(\Omega^{2}_{E}(a\otimes
b)\otimes c)+ \chi_{b,\;c}\Omega^{2}_{E}(\Omega^{2}_{E}(a\otimes c)\otimes b)
+ \chi_{a,\;b+c}\chi^{-1}_{a, E(bc)}\Omega^{2}_{E}(a, \Omega^{2}_{E}(b\otimes c))\\
=\Omega^{3}_{E^{2}}(a\otimes b\otimes c) -E\Omega^{3}_{E}(a\otimes b\otimes
c) -\Omega^{3}_{E}(E(a)\otimes b\otimes c)\\
-\chi_{a,\;b}\chi^{-1}_{a,\;E(b)}\Omega^{3}_{E}(a\otimes E(b)\otimes c)
-\chi_{a+b,\;c}\chi^{-1}_{a+b),\;E(c)}\Omega^{3}_{E}(a\otimes b\otimes E(c)).  \end{multline*} \end{lem}
\prof For the right side the lemma \ref{der2} for
$\Omega^{3}_{E}$ with (\ref{Jacobi})
is used. The corresponding terms we can write in the following way:
\begin{multline*}
\Omega^{3}_{E^{2}}(a\otimes b\otimes c)=
E\Omega^{2}_{E}(a,bc)+\Omega^{2}_{E}(Ea,bc) +
\chi_{a,bc}\chi^{-1}_{a,E(bc)}\Omega^{2}_{E}(a,E(bc))\\
-E(\Omega^{2}_{E}(a,b))c-\Omega^{2}_{E}(Ea,b)c-\chi_{a,b}\chi^{-1}_{a,E(b)}c
-\chi_{b,c}E(\Omega^{2}_{E}(a,c)b)\\-\chi_{b,c}\Omega^{2}_{E}(Ea,c)b
-\chi_{b,c}\chi_{a,c}\chi^{-1}_{a,Ec}\Omega^{2}_{E}(a,Ec)b,
\end{multline*}
and
\begin{multline*}
E\Omega^{3}_{E}(a\otimes b\otimes c)=
E\Omega^{2}_{E}(a,bc)-E(\Omega^{2}_{E}(a,b)c)-
\chi_{b,c}E(\Omega^{2}_{E}(a,c)b),\\
\Omega^{3}_{E}(Ea\otimes b\otimes c)=
\Omega^{2}_{E}(Ea,bc)-\Omega^{2}_{E}(Ea,b)c-\chi_{b,c}\Omega^{2}_{E}(Ea,c)b,\\
\Omega^{3}_{E}(a\otimes Eb\otimes c)=
\Omega^{2}_{E}(a,E(b)c)-\Omega^{2}_{E}(a,E(b))c-
\chi_{E(b),c}\Omega^{2}_{E}(a,c)E(b),\\
\Omega^{3}_{E}(a\otimes b\otimes Ec)=\Omega^{2}_{E}(a,bE(c))-
\Omega^{2}_{E}(a,b)E(c)-\chi_{b,E(c)}\Omega^{2}_{E}(a,E(c))b.
\bbox\end{multline*}
Denote by $D\in\End(A)$ the special linear map such that
\begin{equation}\label{con} \Omega^{3}_{D}=0,\qquad
D^{2}=0\qquad\text{and}\qquad D(1)=0.  \end{equation}

\begin{prop}[Jacobi identity]\label{Jac}  For the linear map (\ref{con})
$D\in\End(A)$ and bicharacters \ref{OmegaE2} , \ref{Jacobi} we have an
identity for the binary operations $$ \Omega^{2}_{D}(\Omega^{2}_{D}(a\otimes
b)\otimes c)+ \chi_{b,\;c}\Omega^{2}_{D}(\Omega^{2}_{D}(a\otimes c)\otimes
b)+ \chi_{a,\;b+c}\chi^{-1}_{a, D(bc)}\Omega^{2}_{D}(a,
\Omega^{2}_{D}(b\otimes c))=0 $$ \end{prop} \prof The proof is based on the
lemma \ref{Jacobi1} with (\ref{con}). $\bbox$\\

From the binary operations $\Omega^{2}$ we can construct a braided
derivations of a multiplication of the algebra $A$ by a term
$d_{a}=\Omega^{2}_{D}(a,\cdot)$.  \begin{prop}\label{endd} For the linear map
(\ref{con})  $D\in\End(A)$ we have \begin{gather*}\label{der}
\Omega^{2}_{D}(a\otimes bc)= \Omega^{2}_{D}(a\otimes
b)c+\chi_{b,\;c}\chi^{-1}_{\Omega^{2}_{D}(a,\; c),\;b}
b\Omega^{2}_{D}(a\otimes c).  \end{gather*} \end{prop} \prof This follows
from the lemma \ref{der2}

Consider a function $c:G\times G\rightarrow \k$, $\forall g_{1},\;g_{2}\in G,
c(g_{1},\;g_{2})\equiv c_{g_{1},\;g_{2}}\in\k$.  For the linear map $D$ the
bracket $[\cdot,\cdot]_{D}$ is defined via the binary operation
$\Omega^{2}_{E}$ \begin{gather}\label{cfun}
c_{|a|,|b|}[a,\;b]_{D}:=\Omega^{2}_{D}(a,\;b).  \end{gather} Let a function
$c$ satisfies \begin{gather}\label{cfun1} c(|a|,|b|)=-c(|D(a)|, |b|),\qquad
c(|a|,|b|)=-c(|a|, |D(b)|)\in\k. \end{gather}
\begin{prop}[Derivation]\label{deriv} For the bracket (\ref{cfun}) with
(\ref{cfun1}) the linear map (\ref{con}) $D\in\End(A)$ is a braided
derivation, $$ D[a,\;b]_{D}=[D(a),\; b]_{D}
+\chi_{a,\;b}{\chi}^{-1}_{a,\;D(b)} [a,\; D(b)]_{D}.  $$ \end{prop} \prof For
the special linear map (\ref{con}) the right side of the lemma \ref{Der1} is
equal to zero.$\bbox$\\

\end{document}